\documentclass[a4paper, 12pt, leqno]{amsart}

\usepackage{amsmath,amsfonts,xcolor}
\usepackage{amsmath}
\usepackage{graphicx}
\usepackage{rotating}
\usepackage{caption}
\usepackage{fancyhdr}
\usepackage{hyperref}
\usepackage{amsaddr}

\begin{document}
\pagestyle{plain}
\bibliographystyle{plain}
\newtheorem{theo}{Theorem}[section]
\newtheorem{lemme}[theo]{Lemma}
\newtheorem{cor}[theo]{Corollary}
\newtheorem{defi}[theo]{Definition}
\newtheorem{prop}[theo]{Proposition}
\newtheorem{problem}[theo]{Problem}
\newtheorem{remarque}[theo]{Remark}
\newcommand{\beq}{\begin{eqnarray}}
\newcommand{\enq}{\end{eqnarray}}
\newcommand{\be}{\begin{eqnarray*}}
\newcommand{\en}{\end{eqnarray*}}
\newcommand{\ben}{\begin{eqnarray*}}
\newcommand{\enn}{\end{eqnarray*}}
\newcommand{\Td}{\mathbb T^d}
\newcommand{\Rd}{\mathbb R^n}
\newcommand{\R}{\mathbb R}
\newcommand{\N}{\mathbb N}
\newcommand{\Sn}{\mathbb S}
\newcommand{\Zd}{\mathbb Z^d}
\newcommand{\Linf}{L^{\infty}}
\newcommand{\dt}{\partial_t}
\newcommand{\Dt}{\frac{d}{dt}}
\newcommand{\Dtt}{\frac{d^2}{dt^2}}
\newcommand{\demi}{\frac{1}{2}}
\newcommand{\vf}{\varphi}
\newcommand{\epu}{_{\varepsilon}}
\newcommand{\ep}{^{\varepsilon}}
\newcommand{\bfi}{{\mathbf \Phi}}
\newcommand{\bpsi}{{\mathbf \Psi}}
\newcommand{\bx}{{\mathbf x}}
\newcommand{\dis}{\displaystyle}
\newcommand{\ds}{\partial_s}
\newcommand{\dss}{\partial_{ss}}
\newcommand{\dx}{\partial_x}
\newcommand{\dxx}{\partial_{xx}}
\newcommand{\dy}{\partial_y}
\newcommand{\dyy}{\partial_{yy}}
\newcommand{\dtt}{\partial_{tt}}
\newcommand {\g}{\`}
\newcommand{\E}{\mathbb E}
\newcommand{\bQ}{\mathbb Q}
\newcommand{\1}{\mathbb I}
\newcommand{\bF}{\mathbb F}
\newcommand{\F}{\cal F}
\newcommand{\bP}{\mathbb P}
\let\cal=\mathcal
\newcommand{\lb}{\langle}
\newcommand{\rb}{\rangle}
\newcommand{\bu}{\bar{u}}
\newcommand{\uu}{\underline{u}}
\newcommand{\bv}{\bar{v}}
\newcommand{\uv}{\underline{v}}
\newcommand{\cS}{{\cal S}}
\newcommand{\cSf}{\cS_\Phi}
\newcommand{\bw}{\bar{w}}
\newcommand{\uw}{\underline{w}}
\newcommand{\bsig}{\bar{\sigma}}
\newcommand{\usig}{\underline{\sigma}}
\newcommand{\mT}{\text {c-exp}}
\newcommand{\mTT}{{\rm c^*}\text{-exp}}

\def\greg#1{{\color{red}#1}}
\def\red#1{{\color{red}#1}}
\def\blue#1{{\color{blue}#1}}
\def\green#1{{\color{green}#1}}


\author{Gr\'egoire Loeper$^{1}$, Neil S Trudinger$^2$} 
\address{$^{1}$Monash University, School of Mathematical Sciences\\
$^2$ The Australian National University, Mathematical Sciences Institute }
\email{gregoire.loeper@monash.edu, neil.trudinger@anu.edu.au}
\title{Weak formulation of the MTW condition and convexity properties of potentials}
\date{\today}

\maketitle

\begin{abstract}
We simplify the geometric interpretation of the weak Ma-Trudinger-Wang condition for regularity in optimal transportation and provide a largely geometric proof of the global c-convexity of locally $c$-convex potentials when the cost function $c$ is only assumed twice differentiable.
\end{abstract}


\section{Introduction}
We consider a cost function $c$ defined on the product $\Omega\times\Omega^*$ of two domains $\Omega, \Omega^*$  in Euclidean space $\R^n$. For a mapping $\phi:\Omega\to \R$ we define its c-transform $\phi^c: \Omega^* \to \R$ by
\ben
\forall y\in \Omega^*, \phi^c(y) = \sup_{x\in \Omega}\{-\phi(x) - c(x,y)\}.
\enn
Conversely we define the $c^*$- transform of $\psi:\Omega^* \to \R$.
A c-convex potential has at every point $x\in \Omega$ a c-support, i.e., there exists $y\in \Omega^*, \psi=\psi(y)\in \R$ such that
\ben
\forall x'\in \Omega, \phi(x') \geq -\psi(y) -c(x',y),
\enn
with equality at $x'=x$. It follows from this definition that 
\[
\phi(x)=\sup_{y\in \Omega^*}\{-\psi(y)-c(x,y)\}
\] and that  $\phi$
can be obtained as the $c^*$ transform of  $\psi:\Omega^*\to \R$. It then turns out that $\psi=\phi^c$. 
For $\phi$ a c-convex potential, and $\phi^c$ its c-transform, we define as in \cite{L6} the contact set as a set valued map $G_\phi$ given by
\ben
G_\phi(y) = \{x: \phi(x)+ \phi^c(y) = -c(x,y)\}.
\enn
for $y \in \Omega^*$. 
We will also use the notions of c-segment, c-convexity of domains. Whenever needed, we will refer to the conditions {\bf A1}, {\bf A2}, {\bf A3}, {\bf A3w} that have been introduced in \cite{MTW, TW}.
One of the main features  of this paper is that we will assume throughout  that the cost function $c$ is globally $C^2(\Omega\times\Omega^*)$, without any further explicit smoothness hypotheses.
As usual we will use subscripts to denote partial derivatives of c with respect to variables $x\in \Omega$ and subscripts preceded by a comma to denote partial derivatives with respect to $y\in\Omega^*$, so that in particular $c_x, c_i, c_{,y}, c_{,j}, c_{i, j}$ denote the partial derivatives of $c$ with respect to $x, x_i, y, y_j, x_i y_j$. We also use $c^{x,y} =[c^{i,j}]$ to denote the inverse of the matrix $c_{x,y} =[c_{i, j}]$. We further assume throughout the paper that $c$  satisfies the assumptions ${\bf A1}, {\bf A2}$ of \cite{MTW}, that is for all $x\in \Omega$ the mapping $y \to -c_x(x,y)$ is injective, that the dual counterpart holds and the matrix $c_{x,y}$ is not singular. We also introduce what will be a weak form of assumption {\bf A3w}:

\begin{defi}\label{A3v}
The cost function satisfies {\bf A3v} if:
for all $x, x_1 \in \Omega$ and $y_0, y_1 \in \Omega^*$, for all $\theta \in (0,1)$, with \[ c_x(x,y_\theta)=\theta  c_x(x,y_1)+(1-\theta) c_x(x,y_0),\]
there holds 
\be
& &\max\{-c(x,y_0) + c(x_0,y_0), -c(x,y_1) + c(x_0, y_1) \} \\
&\geq & -c(x,y_\theta) + c(x_0, y_\theta)  + o(|x-x_0|^2),
\en
where the term $o(|x-x_0|^2)$ may depend on $\theta$.
\end{defi}

From \cite{L6} it is known that when the cost function is $C^4$, {\bf A3v} is equivalent to {\bf A3w}.


Our main result is the following:

\begin{theo}\label{main1}
Let $c : \Omega \times \Omega^* \to \R$ be a $C^2$ cost-function  satisfying  {\bf A1}, {\bf A2} with $\Omega, \Omega^*$ c-convex with respect to each other. Assume that 
\begin{enumerate}
\item[(i)] $c$ satisfies {\bf A3v}.
\end{enumerate}
Then
\begin{enumerate}
\item[(ii)] for all $y_0, y_1 \in \Omega^*$,  $\sigma\in \R$, the set $U =\{x \in\Omega: c(x,y_0)-c(x,y_1)\leq \sigma\}$ is c-convex with respect to $y_0$,
\item[(iii)] for all $\phi$ c-convex, $x\in \Omega$, $y\in \Omega^*$, the contact set $G_\phi(y)$ and its dual $G_{\phi^c}(x)$ are connected,
\item[(iv)] any locally c-convex function in $\Omega$ is globally c-convex.
\end{enumerate}
\end{theo}
\textsc{Remark.} The novelty of the result lies in the way it is obtained; at no point do we have to differentiate the cost function $c$.
Hence the computations from previous proofs \cite{KimMcCann,vi-loep,TW2}, in the case when $c\in C^4$,  do not have to be reproduced. The proof will be based on a purely geometric interpretation of condition {\bf A3v}.


\section{Proof of Theorem \ref{main1}}  
In what follows we will use the term $c-$exponential ($\mT$), as in \cite{L6}, to denote the mapping in condition ${\bf A1}$, that is

\[y=\mT_x(p) \Leftrightarrow -c_x(x,y)=p.\]
We recall also that 
\[D_p( \mT_x)= - c^{x,y}.\]

The core of the proof lies in the following two lemmas,

\begin{lemme}[c-hyperplane lemma]\label{c-hyperplane}
Let $x_0 \in \Omega$, $y_0, y_1\in\Omega^*$ and let $y_\theta = \mT_{x_0} p_{\theta}$ where $p_\theta = (1-\theta)c_x(x_0, y_0)+\theta c_x(x_0, y_1)$, $0\le \theta\le 1$, denote a point on  the c-segment from $y_0$ to $y_1$, with respect to $x_0$. 
Consider, for $\theta > 0$, the section,
\be
S_\theta = S(x_0, y_0, y_\theta): = \{x \in \Omega: c(x,y_0) -c(x_0, y_0) \leq c(x,y_\theta) -c(x_0, y_\theta) \}.
\en
Then as $\theta$ approaches 0, $\partial S_\theta\cap\Omega$ converges to $H_0$,  the $c^*$-hyperplane with respect to $y_0$, passing through $x_0$, with c-normal vector $p_1-p_0$, given by
\[
H_0= H_0( x_0, y_0, y_1) = \{x\in \Omega: -c^{x,y}(x_0,y_0) (p_1-p_0)\cdot [c_{,y}(x,y_0)-c_{,y}(x_0,y_0)] = 0 \}
\]

\end{lemme}
\textsc{Proof.} Locally around $\theta=0$, the equation of $\partial S_\theta$ reads
\ben
[c_{,y}(x,y_0)-c_{,y}(x_0,y_0)] \cdot (y_\theta-  y_0) = o(\theta).
\enn 
Passing to the limit as $\theta$ goes to $0$, we obtain 
\ben 
 [c_{,y}(x,y_0)-c_{,y}(x_0,y_0)] \cdot \partial_\theta y_\theta = 0,
\enn 
which gives the desired result, since 
\ben
\partial_\theta y_\theta &=&  -c^{x,y}(x_0,y_0)(p_1-p_0).
\enn

$\hfill \Box$

\textsc{Remark.} We call $H_0$ a c-hyperplane with respect to $y_0$ because if we express $x$ as $\mTT_{y_0}(q)$ then 
\[H_0=\mTT_{y_0}(\tilde H_0),\]
or equivalently 
\[\tilde H_0 = -c_{,y}(\cdot,y_0)(H_0),\]
where
\[
\tilde H_0= \{q\in c_{,y}(\cdot, y_0)(\Omega): c^{x,y}(x_0,y_0)(p_1-p_0).(q-q_0)=0\}, \quad q_0 = -c_{,y}(x_0,y_0)
\]
Therefore, $H_0$ is the image by $\mTT_{y_0}$ of a hyperplane.
\vspace{1cm}

\textsc{Remark.} We will define in the same way, (replacing $0$ by $\theta$ and $\theta$ by $\theta^\prime$), the section  $S_{\theta, \theta^\prime} $, for 
$\theta^\prime \in (\theta,1)$, and the $c^*$ -hyperplane, $H_\theta = \lim_{\theta^\prime \rightarrow \theta} S_{\theta, \theta^\prime}$.

\vspace{1cm}
The following lemma is then the second main ingredient of the proof: it says that the c-convexity of  $S_\theta$ is non-decreasing with respect to $\theta$; (note that the previous lemma asserts that the ${\rm c}$-convexity of $S_\theta$ vanishes at $\theta=0$).
\begin{lemme}\label{local-convex}
Assume that $c$ satisfies {\bf A3v}.Then 
the second fundamental form of $\partial S_\theta$ at $x_0$ is non-decreasing with respect to $\theta$, for $\theta$ in $(0,1]$.
\end{lemme}

\textsc{Proof.}  Consider
$$h_\theta = c(x,y_0) -c(x,y_\theta) - c(x_0,y_0)  + c(x_0, y_\theta).$$
Note that $h_\theta$ is a defining function for $S_\theta$ in the sense that $S_\theta = \{x\in \Omega: h_\theta\leq 0\}$.

Note also that at $x=x_0$ we have
$h_\theta(x_0)=0$ for all $\theta$ and the set 
\[
\{\partial_x h_\theta |_{x=x_0}, \theta \in [0,1]\}
\] is a line. Therefore all the sets $\partial S_\theta$ contain $x_0$ and have the same unit normal at $x_0$.

Then we note that property {\bf A3v} is equivalent to the following: locally around $x_0$ we have

\beq\label{h-order} 
h_\theta \leq \max\{h_1,0\} + o(|x-x_0|^2).
\enq
(To see this, we just subtract $c(x_0,y_0)-c(x,y_0)$ from both sides of the inequality {\bf A3v}).

Then (\ref{h-order}) implies that the second fundamental form of $\partial S_\theta$ cannot strictly dominate the second fundamental form of
 $\partial S_1$ in any tangential direction at $x_0$.  By changing $y_1$ into $y_{\theta'}$ for $\theta' \geq \theta$, this implies that the second fundamental form of $\partial S_\theta$ is non-decreasing with respect to $\theta$.

$\hfill \Box$

\textsc{Remark.} We remark that analytically the conclusion of Lemma \ref{local-convex} can be expressed as a co-dimension one convexity of the matrix $ A(x, p) = -c_{xx}( x,  \mT_{x_0}(p))$ with respect to $p$, in the sense that the quadratic form $A\xi.\xi$ is convex on line segments in $p$ orthogonal to $\xi$ or more explicitly:
\beq\label{analyticA3}
\Big[A_{ij}(x,p_\theta) - (1-\theta) A_{ij}(x,p_0) - \theta A_{ij}(x,p_1)\Big] \xi_i\xi_j \leq 0,
\enq
for all $\xi \in \R^n$ such that $\xi\cdot(p_1-p_0)=0$, which, for arbitrary $y_0, y_1\in \Omega^*$, is clearly equivalent to  {\bf A3w} when $c\in C^4$.

%
%
%
%
%
%
%
%

\vspace{1cm}

We now deduce assertion $(ii)$ in Theorem 1.2 from ${\bf A3v}$; this will be done in several steps.
\\
Step 1. Uniform boundedness of the section's curvature (including c-hyperplanes)
\\
From the previous corollary,
it follows that $\theta \to c_{xx}(x_0,y_\theta)\xi_i\xi_j$ is convex and therefore Lipschitz, and for a.e. $\theta \in [0,1]$, 
$$A=\partial_\theta c_{x_ix_j}(x_0,\mT_{x_0}(p_{\theta}))\xi_i\xi_j$$ exists and is equal to $\lim_{\theta' \to \theta}B(\theta,\theta')$ where
\ben
B(\theta,\theta')=\frac{\left(c_{x_ix_j}(x_0,\mT_{x_0}(p_{\theta'})-c_{x_ix_j}(x_0,\mT_{x_0}(p_{\theta})\right)\xi_i\xi_j}{\theta'-\theta}.
\enn
The first term $A$ would be  the curvature of $H_\theta$ if it exists. 
The second term $B$ in the limit is the curvature of $S_{\theta, \theta'}$. We can deduce right away that the curvature of $S_{\theta,\theta'}$ remains uniformly bounded at $x_0$ thanks to \eqref{analyticA3}. 
Now this reasoning can be extended to any point $x_1\in \partial S_{\theta, \theta'}$, although the c-segment between $y_\theta$ and $y_\theta'$ will be with respect to $x_1$, but the conclusion that the curvature of $S_{\theta,\theta'}$ at $x_1$ is uniformly bounded remains. Therefore the curvature of all sections is uniformly bounded so as the uniform limit of $\partial S_{\theta, \theta'}$, 
$H_\theta$ is a $C^{1,1}$ hypersurface, and therefore has a curvature a.e. given by $A$.

\vspace{1cm}
 Step 2. Local convexity
Wherever $A$ is well defined, the curvature of $H_\theta$ is equal to $A$.
Moreover, for $\theta'> \theta$, the second fundamental form of $\partial S_{\theta, \theta'}$ dominates a.e. the one of $H_\theta$.

Let us define the hypersurfaces
\[
P_m=\{x\in \Omega, c(x,y_0)-  c(x,y_1)=m\}, m\in \R.
\]
By standard measure theoretical arguments, the previous result implies the following:
\begin{lemme}
For a.e. $y_0, y_1, m$ there holds at ${\cal H}^{n-1}$ every point $x_0$ on $P_m(y_0,y_1)$, that 
\begin{itemize}
\item[-] the second fundamental form (SFF) of $H_0(x_0,y_0,y_1)$ at $x_0$ is well defined, let us call it $A$, equivalently $H_0(x_0,y_0,y_1)$ is twice differentiable (as a hypersurface)
\item[-]  $A$ is dominated by the SFF of $\partial S_1(x_0,y_0)$
\item[-] going back to the tangent space (i.e. composing with $c_{,y}(\cdot, y_0)$), the second fundamental form of $c_{,y}(\cdot, y_0)(\partial S_1(x_0,y_0))$ dominates the null form.   
\end{itemize}
\end{lemme}

We now conclude the local convexity.
Starting from a point $x_0$ where $H_0$ and $S_0$ are tangent.  Both are defined by $x_0,y_0,y_1$.
Representing $\partial S_1$ and $H_0$ as graphs over $\R^{n-1}$, and  we denote by $s_1$ and $s_0$ the corresponding functions.
We assume $x_0=0$, and that both graphs have a flat gradient at $0$. For $x^\prime \in R^{n-1}$ we have
\ben
s_i(x')=|x'|^2\int_0^1\partial_{\nu\nu}s_i(\theta x')(1-\theta)d\theta , \quad i = 0,1, 
\enn
where $\nu$ is the appropriate unit vector. By the definition of $H_0$, at a given point $x^*=(x', h_0(x'))$, $H_0$ is tangent to 
\[
S_1^* = S(x^*,y_0, y_1^*), \quad y_1^* = \mT[x^*,-c_x(x^*,y_0)+ c_{x,y}(x^*,y_0)c^{x,y}(x_0,y_0)(p_1-p_0)].
\]
For almost every choice of $x_0$ there will hold  for a.e. $x'$ that 
\ben
\partial_{\nu\nu}s_0(x') &\leq& \partial_{\nu\nu} s^*_1(x')\\
&\leq &  \partial_{\nu\nu} s_1(x') + \varepsilon(x'-0)
\enn
with $\lim_0\varepsilon=0$, depending on the continuity of $c_{xx}, c_{x,y}$.
Therefore
\ben
s_0(x') &\leq& |x'|^2 (\int_0^1  \partial_{\nu\nu}s_1(\theta x')(1-\theta)d\theta + \varepsilon(x'))\\
&\leq& s_1(x') + \varepsilon(x'-x_0)|x'|^2.
\enn
Going now in the tangent space, for $q'$ in a well chosen $n-1$ subspace, and $\pi$ the projection on $\{x_n=0\}$, we call $x(q') = \pi(\mTT(y_0,q'))$ and we have  
\ben
s_0(x(q')) \leq s_1(x(q')) + \varepsilon(x(q')-x_0)|x(q')|^2,
\enn
$s_0(x(q'))$ is an affine function, $s_1(x(q'))$ defines the image of $\partial S_1$ by $c_{,y}$ and $\varepsilon(x(q'))|x(q')|^2 \leq \tilde\varepsilon(q'))|q'|^2$ for some $\varepsilon'$. 
For a.e. choice of $x_0$, this holds for a.e $q'$. More importantly the $\varepsilon'$ is (locally) uniform.
This implies the convexity through the following lemma
\begin{lemme}
Let $s$ be $C^1$. Assume that for some continous $\varepsilon(\cdot)$ with $\varepsilon(0)=0$, there holds for almost every $x_0,x$
\[
s(x)\geq l_{x_0}(x) - \varepsilon(x-x_0)|x-x_0|^2
\]
$l_{x_0}$ being the tangent function  at $x_0$, then $s$ is convex.
\end{lemme}
{\it Proof.} Elementary, both sides of the inequality are continuous in $x,x_0$, so this holds in fact everywhere.

$\hfill \Box$

\textsc{Remark.} For a proof of local convexity without using Lemma 2.3 the reader is referred to \cite{ LT2}.


\vspace{1cm} Global convexity
 To complete the proof of assertion (ii), we need to  show that the set $\tilde S_1$ is connected. The proof goes as follows, and it is very close to the argument of \cite{TW}, Section 2.5. Let $\sigma$ be a constant, and assuming that the set
\[ \{c(x,y_0) -c (x,y_1)\leq \sigma\}\] has two disjoint components, we let $\sigma$ increase until the two components touch in a $C^1$ c -convex subdomain $\Omega' \subset \subset \Omega$. 
From the local convexity property this can only happen on the boundary of $\Omega'$. 
At this point, say $x_1$ there holds locally that \[ c(x,y_0)-c(x,y_1) \leq \sigma\] on $\partial\Omega'$ and for $x\ep=x_1- \varepsilon \nu$, $\nu$ the outer unit normal to $\Omega'$, 
\[ c(x\ep,y_0)-c(x\ep,y_1) > h.\]
 This implies that
\[ c(x,y_0)-c(x,y_1) \geq h\] 
is locally c-convex around $x_1$, a contradiction, and from this we deduce that $S_1$ can have at most one component. Since a connected locally convex set in Euclidean space must be globally convex, we thus deduce that $S_1$ is globally c-convex.
 
 
 $\hfill \Box$

\subsection{An analytical proof for a smooth cost function} 
If a $C^2$ domain $\Omega$ is defined locally by $\varphi  > 0$, its
local c-convexity with respect to $y_0$, for $c \in C^3$, is expressed by
\be
\Big[\varphi_{ij} + c_{ij,k}c^{k,l}(\cdot, y_0)\partial_l\varphi\Big] \tau_i\tau_j \geq 0,
\en   
 or equivalently
\be
\Big[\varphi_{ij} + \partial_pA^{ij}.\partial\varphi\Big] \tau_i\tau_j \geq 0
 \en 
 for all $\tau \in \partial\Omega$ \cite{TW2}.
Plugging $\varphi(x) = c(x,y_0) - c(x,y_1) - h$ into this inequality, we obtain immediately from (2) that $S_1$ is locally $c$-convex with respect to $y_0$.
More generally this argument proves Theorem 1.2 when we assume additionally that the form  $A\xi.\xi$ is differentiable with respect to $p$ in directions orthogonal to $\xi$.





\subsection{Connectedness of the contact set}

This new characterization implies right away the $c$-convexity  of the global $c$-sub-differential, ($c$-normal mapping).  We prove now that $(i)$ implies $(iii)$.
 
For $\phi$ c-convex, we have 
\beq\label{base}
\phi(x) &=& \sup_y \{-\phi^c(y)-c(x,y)\},\\
\phi^c(y) &=& \sup_x \{-\phi(x)-c(x,y)\}.
\enq
Then
\be
\{\phi(x) \leq -c(x,y_0) + h \}&=& \cap_y \{x: -\phi^c(y)-c(x,y) \leq -c(x,y_0) + h\}\\
&=&  \cap_y \{x: c(x,y_0) \leq c(x,y) - h + \phi^c(y)\}.
\en
Therefore $\{\phi(x) \leq -c(x,y_0) + h\}$ is an intersection of c-convex sets and hence also c-convex.
We then have 
\be
G_\phi(y) &=& \{x, \phi(x)=-c(x,y) -\phi^c(y)\}\\
&=& \{x, \phi(x)\leq -c(x,y) -\phi^c(y)\},
\en
and hence $G_\phi(y)$ is a c-convex set. To show the dual conclusion, we may rewrite assertion (ii) as:
for all $y, y_1 \in \Omega^*$, $x_0, x_1 \in \Omega^*$ and $\theta \in (0,1)$, with 
\[ c_y(x_\theta,y)=\theta  c_y(x_1,y)+(1-\theta) c_y(x_0,y),\]
there holds 
\be
& &\max\{-c(x_0,y) + c(x_0,y_0), -c(x_1,y) + c(x_1, y_0) \} \\
&\geq & -c(x_\theta,y) + c(x_\theta, y_0).
\en
Since this shows in particular that {\bf A3v} is invariant under duality we complete the proof of assertion (iii). Moreover as a byproduct of this argument we also see that the sets $S_\theta$ are non-increasing with respect to $\theta$ and that {\bf A3v} holds without the term $o(|x-x_0|^2)$.

$\hfill \Box$

\subsection{Local implies global} 
We prove that $(ii)$ implies $(iv)$. 
We consider $\phi$ a locally c-convex function, i.e, $\phi$ has at every point a local c-support. Locally, $\phi$ can be expressed as
\ben
\phi(x) = \sup_{y\in \omega}\{-\psi(y)-c(x,y)\},
\enn
for some $\omega(x) \subset \Omega^*$
(if $\phi$ was globally c-convex there would hold that $\omega \equiv \Omega^*$ and $\psi$ would be equal to $\phi^c$ ). It follows that 
the level sets 
\[
S_{m,y_0}=\{x: \phi(x)+ c(x,y_0) \leq m\}
\] are locally c-convex with respect to $y_0$ for any $y_0$.  We obtain that $-\partial_yc(S_{m,y_0}, y_0)$ is locally convex. 
Reasoning again as in the proof of the global convexity in point $(ii)$ (i.e. increasing $m$ until two components touch), we obtain that, for $\phi$ locally c-convex,  $-\partial_yc(S_{m,y_0}, y_0)$ is globally convex for all $y_0$.
This implies in turn the global c-convexity of $\phi$, following Proposition 2.12 of \cite{L6}. As already mentioned, this part is very similar to the argument of \cite{TW2}, section 2.5.

Finally we remark that the arguments in this paper extend to generating functions as introduced in \cite {T2014} and also provide as a byproduct an alternative geometric proof of the invariance of condition {\bf A3w} under duality  to the more complicated calculation in \cite{T2014}. The resultant convexity theory is presented in \cite{LT2}.

$\hfill\Box$

\bibliography{../biblio-greg} 

\def\cprime{$'$} \def\cprime{$'$}
\begin{thebibliography}{1}

\bibitem{KimMcCann}
Y.-H. Kim and R.~J. McCann.
\newblock Continuity, curvature, and the general covariance of optimal
  transportation.
\newblock {\em J. Eur. Math. Soc. (JEMS)}, 12:1009--1040, 2010.

\bibitem{L6}
G.~Loeper.
\newblock On the regularity of solutions of optimal transportation problems.
\newblock {\em Acta Mathematica}, 202(2):241--283, 2009.

\bibitem{LT2}
G.~Loeper and N.S. Trudinger.
\newblock On the convexity theory of generating functions.
\newblock {\em in preparation}.

\bibitem{vi-loep}
G.~Loeper and C.~Villani.
\newblock Regularity of optimal transport in curved geometry: the non-focal
  case.
\newblock {\em Duke Math. Journal}, 151(3):431--485, 2010.

\bibitem{MTW}
X.-N. Ma, N.S. Trudinger, and X.-J. Wang.
\newblock Regularity of potential functions of the optimal transport problem.
\newblock {\em Arch. Ration. Mech. Anal.}, 177(2):151--183, 2005.

\bibitem{T2014}
N.S. Trudinger.
\newblock On the local theory of prescribed {J}acobian equations.
\newblock {\em Discrete Contin. Dyn. Syst.}, 34(4):1663--1681, 2014.

\bibitem{TW2}
N.S. Trudinger and X.-J. Wang.
\newblock On strict convexity and continuous differentiability of potential
  functions in optimal transportation.
\newblock {\em Arch. Ration. Mech. Anal.}, 192(3):403--418, 2009.

\bibitem{TW}
N.S. Trudinger and X.-J. Wang.
\newblock On the second boundary value problem for {M}onge-{A}mp\`ere type
  equations and optimal transportation.
\newblock {\em Ann. Scuola Norm. Sup. Pisa Cl. Sci.}, 5(8):143--174, 2009.

\end{thebibliography}
\vspace{1cm}
\noindent Gregoire Loeper\\
Monash University\\
School of Mathematics\\
9 Rainforest Walk\\
3800 CLAYTON VIC, AUSTRALIA\\
\\
email: gregoire.loeper@monash.edu

\end{document}